\documentclass[a4paper,11pt]{amsart}

\usepackage{amsfonts} 
\usepackage{amsmath} 
\usepackage{amsthm} 

\title[Pointwise Universal Gysin formul\ae{}]{Pointwise Universal Gysin formul\ae{} and Applications towards Griffiths' conjecture}
\author{Simone Diverio}
\address{Simone Diverio \\ Dipartimento di Matematica \lq\lq Guido Castelnuovo\rq\rq{} \\ SAPIENZA Universit\`a di Roma \\ Piazzale Aldo Moro 5 \\ I-00185 Roma.}
\email{diverio@mat.uniroma1.it}
\author{Filippo Fagioli}
\address{Filippo Fagioli \\ Dipartimento di Matematica \lq\lq Guido Castelnuovo\rq\rq{} \\ SAPIENZA Universit\`a di Roma \\ Piazzale Aldo Moro 5 \\ I-00185 Roma.}
\email{fagioli@mat.uniroma1.it} 
\curraddr{Dipartimento di Matematica e Informatica \lq Ulisse Dini\rq{} \\Universit\`a degli studi Firenze \\ Viale Morgagni 67/A  \\ I-50134 Firenze}
\email{filippo.fagioli@unifi.it}
\thanks{The first-named author is partially supported by the ANR Programme: D\'efi de tous les savoirs (DS10) 2015, \lq\lq GRACK\rq\rq, Project ID: ANR-15-CE40-0003ANR, and by the ANR Programme: D\'efi de tous les savoirs (DS10) 2016, \lq\lq FOLIAGE\rq\rq, Project ID: ANR-16-CE40-0008}
\date{\today}

\theoremstyle{plain}
\newtheorem{thm}{Theorem}[section]

\newtheorem*{mainthm*}{Main Theorem}
\newtheorem*{mainapp*}{Main Application}
\newtheorem{lem}[thm]{Lemma}
\newtheorem{prop}[thm]{Proposition}
\newtheorem{quest}[thm]{Question}

\theoremstyle{remark}
\newtheorem{rem}[thm]{Remark}
\newtheorem{ex}[thm]{Example}

\theoremstyle{definition}

\newcommand{\Z}{\mathbb{Z}}
\newcommand{\Q}{\mathbb{Q}}
\newcommand{\R}{\mathbb{R}}
\newcommand{\C}{\mathbb{C}}
\newcommand{\PP}{\mathbb{P}}
\newcommand{\GG}{\mathbb{G}}
\newcommand{\FF}{\mathbb{F}_{\rho}}
\newcommand{\OO}{\mathcal{O}}
\newcommand{\ff}{\mathbf{f}}
\newcommand{\pis}{\pi_{*}}
\newcommand{\End}{\text{End}}
\newcommand{\Span}{\operatorname{Span}}
\newcommand{\Tr}{\text{Tr}}
\newcommand{\Qa}{Q_{\rho}^{\mathbf{a}}}
\newcommand{\Xia}{\Xi_{\rho}^{\mathbf{a}}}

\begin{document}
\bibliographystyle{amsalpha}
\maketitle
	
	\begin{abstract}
	Let $X$ be a complex manifold, $(E,h)\to X$ be a rank $r$ holomorphic Hermitian vector bundle, and $\rho$ be a sequence of dimensions $0 = \rho_0 < \rho_1 < \cdots < \rho_m = r$. Let $Q_{\rho,j}$, $j=1,\dots,m$, be the tautological line bundles over the (possibly incomplete) flag bundle $\mathbb{F}_{\rho}(E) \to X$ associated to $\rho$, endowed with the natural metrics induced by that of $E$, with Chern curvatures $\Xi_{\rho,j}$. We show that the universal Gysin formula \textsl{\`{a} la} Darondeau--Pragacz for the push-forward of a homogeneous polynomial in the Chern classes of the $Q_{\rho,j}$'s also holds pointwise at the level of the Chern forms $\Xi_{\rho,j}$ in this Hermitianized situation.
	
	As an application, we show the strong positivity of several polynomials in the Chern forms of a Griffiths (semi)positive vector bundle not previously known, thus giving some new evidences towards a conjecture by Griffiths, which in turn can be seen as a pointwise Hermitianized version of the Fulton--Lazarsfeld theorem on numerically positive polynomials for ample vector bundles.
	\end{abstract}

	\section{Introduction}
Let $E\to X$ be a rank $r\ge 2$ holomorphic vector bundle over a complex manifold of dimension $n$. Once a sequence $\rho$ of dimensions $0 = \rho_0 < \rho_1 < \cdots < \rho_m = r$ has been fixed, one can consider the (incomplete if $m<r$, or complete when $m=r$) flag bundle $\pi\colon\FF(E)\to X$, which is naturally endowed with $m+1$ tautological rank $\rho_j$ vector bundles $U_{\rho,j}\to\FF(E)$, $j=0,\dots,m$. Out of this, one can form tautological line bundles
$$
Q_{\rho,j}:=\det(U_{\rho,m-j+1}/U_{\rho,m-j}),\quad j=1,\dots,m,
$$
and consider the corresponding cohomology classes $c_1(Q_{\rho,j})$ in the cohomology group $H^{2}\bigl(\FF(E)\bigr)$. 

Given a homogeneous polynomial $F$ in $m$ variables of degree $d_\rho+k$, where $d_\rho$ is the relative dimension of the proper holomorphic submersion $\pi\colon\FF(E)\to X$ and $0\le k\le n$, the proper push-forward 
$$
\pi_*F\bigl(c_1(Q_{\rho,1}),\dots,c_1(Q_{\rho,m})\bigr)
$$ 
gives a cohomology class in $H^{2k}(X)$ which of course needs to be a characteristic class for $E$. It is then a natural issue (which has been considered and settled by several authors in different degrees of generality) to try to determine a closed formula to express this class more or less explicitly as a polynomial, call it $\Phi$, in the Chern or Segre classes of $E$; see \cite{Dam73,Ilo78,KT15,DP17}, just to cite a few. We shall consider in particular here the incarnation of such a formula given by Darondeau--Pragacz in \cite[Theorem 1.1, Proposition 1.2]{DP17}

Now, alongside the cohomological situation, one can ask the analogue in the Hermitian setting as follows. Suppose $E$ is moreover endowed with a smooth Hermitian metric $h$. Then, all the tautological bundles $U_{\rho,j}$ considered before inherit, being subbundles of $\pi^*E$, a Hermitian metric, and so do the determinants of the successive quotients $Q_{\rho,j}$. Thus, the classes $c_1(Q_{\rho,j})$ now have special representatives $\Xi_{\rho,j}$ given by the Chern curvatures of their induced Hermitian metrics. 

Given the homogeneous polynomial $F$ above, one can formally compute it using the $\Xi_{\rho,j}$'s as variables to get a closed $(d_\rho+k,d_\rho+k)$-form on $\FF(E)$, which can be pushed-forward on $X$ \textsl{via} integration along the fibers to obtain a $(k,k)$-form $\pi_*F(\Xi_{\rho,1},\dots,\Xi_{\rho,m})$ on $X$. Such a form is of course a special representative for $\pi_*F\bigl(c_1(Q_{\rho,1}),\dots,c_1(Q_{\rho,m})\bigr)=\Phi\bigl(c_\bullet(E)\bigr)$. Certainly, one can also compute the Chern curvature of $(E,h)$, and by the Chern--Weil theory represent the Chern classes of $E$ by the Chern forms $c_j(E,h)$ associated to its Chern curvature. In this way, $\Phi\bigl(c_\bullet(E,h)\bigr)$ is a special representative for $\Phi\bigl(c_\bullet(E)\bigr)$, too. Therefore, \textsl{a priori} $\pi_*F(\Xi_{\rho,1},\dots,\Xi_{\rho,m})$ and $\Phi\bigl(c_\bullet(E,h)\bigr)$ differ by an error term which is an exact $2k$-form.

Our main result (see Theorem \ref{thm: darondeau-pragacz for differential forms}) can be now summarized by saying that

\begin{mainthm*} We have the equality
$$
\pi_*F(\Xi_{\rho,1},\dots,\Xi_{\rho,m})=\Phi\bigl(c_\bullet(E,h)\bigr).
$$
\end{mainthm*}

So, in fact, there is no error term at all: this generalizes previous result of \cite{Mou04,Gul12,Div16} which concerned the case of projectivized bundles, corresponding here to the special weight $\rho=(0,1,r)$ (or, dually, $\rho=(0,r-1,r)$). In other words, the universal Gysin formul\ae{} to compute the push-forwards in cohomology from the flag bundle can be used \textsl{verbatim} to compute pointwise the push-forwards for differential form constructed from the Chern--Weil theory in the Hermitian situation.

\medskip

The second part of the paper is devoted to an application of our Main Theorem to a positivity issue, in the same spirit of \cite{Gul12}, as follows. Suppose $(E,h)\to X$ is a Griffiths (semi)positive vector bundle (for precise definitions, see Section \ref{sec:positivity}). In the seminal paper \cite{Gri69} it is raised the problematic of determine which characteristic forms built from the Chern curvature of $(E,h)$ are positive, and some partial result is given. 

More precisely, it is asked (and expected) whether the Schur forms (and hence, their positive linear combinations) are positive. Griffiths' question was thus a differential, pointwise forerunner of the Fulton--Lazarsfeld theorem \cite{FL83}, which is a cohomological global statement, that characterizes precisely all numerically positive polynomials for ample vector bundles. The Fulton--Lazarsfeld theorem indeed characterizes them exactly as the positive linear combination of Schur polynomials. 

Observe that a Griffiths positive (resp. semipositive) vector bundle on a compact complex manifold is ample\footnote{The converse is not known, but expected to be true.} (resp. nef) and that a cohomology class which can be represented by a positive differential form is numerically positive. Thus, an answer in the affirmative to Griffiths' question would give a stronger Fulton--Lazarsfeld-type statement under the stronger (but conjecturally equivalent) hypothesis of positivity in the sense of Griffiths.

Up to now, very little is known about this question beside the trivial case of (any power of) the first Chern form: other Schur forms known to be positive are the second Chern form \cite{Gri69} (see also \cite{Fag21}), and the signed Segre forms \cite{Gul12} (see again Section \ref{sec:positivity} for a more exhaustive list of related results in literature).

Here, as a consequence of our Main Theorem, we are able to establish the positivity of several new (positive linear combinations of) Schur forms, cf. Section \ref{sec:positivity} and in particular Subsection \ref{subsec:examples}. Namely, thanks to the general curvature formul\ae{} obtained by Demailly in \cite{Dem88} coupled with our universal pointwise Gysin push-forward formula, we obtain the following (see Theorem \ref{thm: positivity of push-forwards by flag bundles}).

\begin{mainapp*}
With notations as above, let $(E,h)\to X$ be a Griffiths semipositive vector bundle. Given a weight $\mathbf{a} \in \Z^m$ such that  $ a_1 \ge \cdots \ge a_m \ge 0 $, we have that the characteristic forms
$$
\pi_*(a_1\,\Xi_{\rho,1}+\cdots+a_m\,\Xi_{\rho,m})^{d_\rho+k}
$$
are strongly positive $(k,k)$-forms on $X$, and moreover positive linear combinations of Schur forms of $(E,h)$.
\end{mainapp*}

This confirms Griffiths' conjecture for those positive linear combinations of Schur forms of $(E,h)$ which can be obtained as wedge products of push-forwards of type $\pi_*(a_1\,\Xi_{\rho,1}+\cdots+a_m\,\Xi_{\rho,m})^{d_\rho+k}$ (see Subsection \ref{subsec:examples} for several concrete examples).

\begin{rem}\label{rem:finski}
We discovered, right before uploading the first version of this article on the arXiv, that a few hours before we completed the redaction of that version it appeared on the arXiv the paper \cite{Fin20} by S. Finski, who proves ---independently of us--- similar results, some of which with techniques not so far from ours. 

In particular, he can prove the positivity of \emph{all} Schur forms but \emph{under the stronger hypothesis} of (dual) Nakano positivity for $(E,h)$ \cite[Theorem 1.1]{Fin20}. He also makes an interesting connection between some open question for the so-called positive semidefinite linear preservers and the original question of Griffiths, showing that they are in fact equivalent problems \cite[Theorem 1.3]{Fin20}. 

Finally, let us remark that \cite[Theorem 7.2]{Fin20} is a special case of our Main Theorem but for complete flag bundles and where, in our notations, $F$ is taken to be just a monomial with some specific decreasing degrees.
\end{rem}

	\section{Flag bundles and their curvature}
	Let $ X $ be a complex manifold of dimension $ n $ and let $ E \to X $ be a holomorphic vector bundle of rank $ r $.
	Fixed a sequence of integers $ \rho=(\rho_0,\ldots,\rho_m) $ of the form $ 0 = \rho_0 < \cdots < \rho_j < \cdots < \rho_m = r $, the \emph{flag bundle} of $ E $ associated to $ \rho $ is the holomorphic fiber bundle
	\[
	\pi \colon \FF(E) \to X
	\]
	where the fiber over $ x \in X $ is the flag manifold $ \FF(E_x) $ of flags
	\[
	\{ 0_{x} \} = V_{x,\rho_0} \subset \cdots \subset V_{x,\rho_j} \subset \cdots \subset V_{x,\rho_m} = E_{x} , \quad \dim_{\C} V_{x,\rho_j} = \rho_j .
	\]
	
	Over $ \FF(E) $ we have a tautological flag
	\begin{equation}\label{eq: tautological filtration vector bundles over flag bundle}
	U_{\rho,0} \subset \cdots \subset U_{\rho,j} \subset \cdots \subset U_{\rho,m}
	\end{equation}
	 of vector subbundles of $ \pi^{*}E $, where the fiber of $ U_{\rho,j} $ over the flag
	\[
	\{ 0_{x} \} \subset \cdots \subset V_{x,\rho_j} \subset \cdots \subset E_{x}
	\]
	is $ V_{x,\rho_j} $; hence $ U_{\rho,j} $ has rank $ \rho_j $.
	The tautological filtration \eqref{eq: tautological filtration vector bundles over flag bundle} allows us to define natural line bundles over $ \FF(E) $ as follows:
	 for $ 1 \le j \le m $ set
	\begin{equation*}
	Q_{\rho,j}:=\det(U_{\rho,m-j+1}/U_{\rho,m-j}).
	\end{equation*}
	For any multi-index $ \mathbf{a} = (a_1,\ldots,a_r) \in \Z^r $ satisfying
	\begin{equation}\label{eq:condition on a}
	a_{r - \rho_{m-j+1} + 1} = a_{r - \rho_{m-j+1} + 2} = \cdots = a_{r - \rho_{m-j} }, \ 1 \le j \le m ,
	\end{equation}
	set, for $0\le j\le m$, $ s_j := r - \rho_{m-j} $ and define
	\[
	\Qa := { Q_{\rho,1} }^{\otimes a_{s_1}} \otimes \cdots \otimes { Q_{\rho,m} }^{\otimes a_{s_m}} .
	\]
	
	In the particular case of complete flag bundles, \textsl{i.e.} $m=r$, we shall drop the subscript $\rho$ and simply write $\mathbb F(E)$, $ U_{j} $, $ Q_{j} $ and $ Q^{\mathbf{a}} $.
	
	\begin{ex}[To keep in mind as a toy case]\label{ex:toycase}
		Suppose that $ E $ is a rank $ r $ vector bundle over $ X $, and consider the flag bundle $ \mathbb{F}_\rho(E) $ corresponding to the sequence $\rho=(0,1,r)$.		
		Then, $ \mathbb{F}_\rho(E) $ coincides with the projectivized bundle of lines $ \PP(E) $ on $ X $.
		By definition, $ U_{\rho,1} $ equals the tautological line bundle $ \OO_{E}(-1) $, while $U_{\rho,2}=\pi^*E$. We thus have the short exact sequence
		$$
		0\to\underbrace{\mathcal O_E(-1)}_{U_{\rho,1}}\to\underbrace{\pi^* E}_{U_{\rho,2}}\to \pi^*E/\mathcal O_E(-1)\to 0,
		$$
		which gives, by taking the determinant, that $ Q_{\rho,1} = \mathcal O_E(1)\otimes\pi^*\det E $ and $ Q_{\rho,2} = \OO_{E}(-1) $.
		Therefore,
		$$
		Q_{\rho}^{(a_1,a_2)}=\mathcal O_E(a_1-a_2)\otimes\bigl(\pi^*\det E\bigr)^{\otimes a_1}.
		$$
	\end{ex}
	
	Suppose that on $ E $ is given a Hermitian metric $ h $, then we have an induced pull-back metric on $ \pi^{*}E $ which endows all the vector bundles $ U_{\rho,j} $ with the restriction metric and all the line bundles $ Q_{\rho,j} $ with the determinant of the quotient metric.
	Consequently, we have natural metrics induced on all the line bundles $ \Qa $.
	In order to simplify the notation, by a slight abuse, denote by $ h $ all these mentioned metrics.
	
	Fix a point $ x_0 \in X $, local holomorphic coordinates $ (z_1,\ldots,z_n) $ on an open set  of $X $ centered at $ x_0 $ and a point $ (x_0,\ff_0) \in \FF(E) $.
	We can always choose a local normal frame $ (e_1,\ldots,e_r) $ of $ E $ at $ x_0 $ such that $ \ff_0 \in \FF(E_{x_0}) $ coincides with the flag
	\[
	\{ 0_{x_0} \} \subset Z_{\rho_1} \subset Z_{\rho_2} \subset \cdots \subset Z_{\rho_{m-1}} \subset E_{x_0} ,
	\]
	where, for $ \ 1 \le j \le m $,
	\[
	Z_{\rho_j} = \Span\bigl\{ e_{r-\rho_j+1}(x_0),e_{r-\rho_j+2}(x_0),\ldots,e_{r}(x_0) \bigr\} .
	\]
	For $z$ in the coordinate open set considered, the basis $\bigl(e_1(z),\ldots,e_r(z)\bigr)$ gives affine coordinates $ \zeta = (\zeta_{\lambda\mu}) $ on the fiber $ \FF(E_z) $, where $ 1 \le \lambda < \mu \le r $ are such that there is an integer $ \ell=1,\ldots,m-1 $ with $ \lambda \le s_\ell < \mu $. Such coordinates parameterize flags of the form
	$$
	\{ 0_{z} \} \subset \cdots \subset\Span\bigl\{\epsilon_{r-\rho_j+1}(z,\zeta),\epsilon_{r-\rho_j+2}(z,\zeta),\ldots,\epsilon_{r}(z,\zeta) \bigr\} \subset \cdots  \subset E_{z},
	$$
	where, for $1\le k\le r$,
	$$
	\epsilon_k(z,\zeta)=e_k(z)+\sum\zeta_{\lambda k}\, e_\lambda(z),
	$$
	and the summation is taken over all $1\le \lambda<k$ such that as before there is an integer $\ell=1,\ldots,m-1$ with $ \lambda \le s_\ell <k $.
	
	Summing up, we have constructed in this way local holomorphic coordinates $ (z,\zeta) = \left( z_1,\ldots,z_n, \zeta_{\lambda\mu} \right)$
	on $ \FF(E) $ around $ (x_0,\ff_0) $.
	
	Let $ \Theta(E,h) \in \mathcal{A}^{1,1}(X,\End(E)) $ be the Chern curvature tensor of $ (E,h) $.
	With respect to the local coordinates introduced before, at $ x_0 $ we have
	\begin{align*}
	\Theta(E,h)_{x_0} &= \sum_{\alpha,\beta=1}^{r}\sum_{p,q=1}^{n} c_{pq\alpha\beta}(x_0)\,dz_p \wedge d\bar{z}_q \otimes e_{\alpha}^{\vee} \otimes e_{\beta} \\
	&= \sum_{\alpha,\beta=1}^{r} \Theta_{\beta\alpha}(x_0) \otimes e_{\alpha}^{\vee} \otimes e_{\beta},
	\end{align*}
	where the $ (1,1) $-forms
	\[
	\Theta_{\beta\alpha} = \sum_{p,q=1}^{n} c_{pq\alpha\beta} \,dz_p \wedge d\bar{z}_q
	\]
	are the entries of the curvature matrix.

For a given multi-index $\mathbf a\in\mathbb Z^r$ satisfying condition (\ref{eq:condition on a}), in \cite[Formula (4.9)]{Dem88} the Chern curvature $ \Theta(\Qa,h) $ of $ (\Qa,h) $ at the point $ (x_0,\ff_0) $ is computed and reads
	\begin{equation}\label{eq: curvature q^a pointwise}
	\Theta(\Qa,h)_{(x_0,\ff_0)} = \sum_{\lambda} a_{\lambda} \Theta_{\lambda\lambda}(x_0) + \sum_{\lambda,\mu} (a_{\lambda} - a_{\mu}) d\zeta_{\lambda\mu} \wedge d\bar{\zeta}_{\lambda\mu},
	\end{equation}
	where $ \lambda $ and $ \mu $ in the second summation are chosen as before in order to have that the $\zeta_{\lambda\mu}$'s are defined. 
	
	\begin{rem}\label{rem:decreasing}
	In the particular case where the multi-index is non increasing, which is indeed the case actually considered in \cite[Formula (4.9)]{Dem88}, and which shall be of special interest later, the condition on $\lambda$ and $\mu$ in the second summation above is equivalent to require $a_\lambda>a_\mu$, so that the curvature formula becomes
	\begin{equation}\label{eq: curvature q^a pointwise decreasing}
	\Theta(\Qa,h)_{(x_0,\ff_0)} = \sum_{\lambda} a_{\lambda} \Theta_{\lambda\lambda}(x_0) + \sum_{a_\lambda>a_\mu} (a_{\lambda} - a_{\mu}) d\zeta_{\lambda\mu} \wedge d\bar{\zeta}_{\lambda\mu}.
	\end{equation}
	Observe that $\sum_{a_\lambda>a_\mu} (a_{\lambda} - a_{\mu}) d\zeta_{\lambda\mu} \wedge d\bar{\zeta}_{\lambda\mu}$ gives a positive definite block for the curvature of $(\Qa,h)$ if and only if the non increasing multi-index $\mathbf a$ is strictly decreasing at each place where it is allowed by condition (\ref{eq:condition on a}), \textsl{i.e.} 
	\begin{equation}\label{eq:star}
	a_{s_1}>a_{s_2}>\cdots>a_{s_m}\tag{$\star$}.
	\end{equation}
	\end{rem}
	
	For $ 1 \le j \le m $, denote by $\mathbf{1}_j $ the element of $ \Z^r $ with $1$'s in the places $ s_{j-1} + 1, s_{j-1} + 2, \ldots, s_{j} $, and $0$'s elsewhere.
	By definition, we have that $ Q_{\rho}^{\mathbf{1}_j} = Q_{\rho,j} $, and thanks to Formula \eqref{eq: curvature q^a pointwise} we can recover the curvature of the line bundle $ Q_{\rho,j} $ at the point $ (x_0,\ff_0) \in \FF(E) $:
	\begin{multline}\label{eq: curvature q_j pointwise}
	\Theta(Q_{\rho,j},h)_{(x_0,\ff_0)} = \sum_{\lambda=s_{j-1} + 1}^{s_j} \Theta_{\lambda\lambda}(x_0) \\ - \sum_{\substack{\lambda=1,\dots,s_{j-1}\\\mu=s_{j-1}+1,\dots,s_j}} d\zeta_{\lambda\mu} \wedge d\bar{\zeta}_{\lambda\mu} + \sum_{\substack{\lambda=s_{j-1}+1,\dots,s_j\\\mu=s_{j}+1,\dots,r}} d\zeta_{\lambda\mu} \wedge d\bar{\zeta}_{\lambda\mu}.
	\end{multline}
	Finally, let
	\[
	\Xi_{\rho,j} := c_1(Q_{\rho,j},h) = \frac{i}{2\pi} \Theta(Q_{\rho,j},h) \in \mathcal{A}_{\R}^{1,1}(\FF(E))
	\]
	be the \emph{first Chern form} of $ (Q_{\rho,j},h) $ which represents the first Chern class $ c_1(Q_{\rho,j}) $, and similarly
	$$
	\Xia := \frac i{2\pi}\Theta(\Qa,h).
	$$
	
	\subsection{Splitting of the tangent bundle and intrinsic expression of the curvature} 
	Now, we rewrite Formula \eqref{eq: curvature q_j pointwise} more intrinsically.
	In order to do this, consider the short exact sequence
	\begin{equation*}
	0 \to \ker(d\pi) \hookrightarrow T_{\FF(E)} \xrightarrow{d\pi} \pi^{*}T_X \to 0
	\end{equation*}
	induced by the differential of $ \pi \colon \FF(E) \to X $, where $ T_{\FF(E)} $ and $ T_X $ are the tangent bundles of $ \FF(E) $ and of $ X $ respectively.
	Recall that $ \ker(d\pi) $ is the \emph{relative tangent bundle}, usually denoted by $ T_{\FF(E)/X} $.
	
	We now define a natural orthogonal splitting of $T_{\FF(E)}$ into a vertical and a horizontal part. In order to do this, observe that for any given weight $\mathbf a$ as described at the end of Remark \ref{rem:decreasing}, by Formula (\ref{eq: curvature q^a pointwise decreasing}) the line bundle $\Qa$ is relatively positive. Hence, $\Xia$ gives a positive definite Hermitian form whenever restricted to the relative tangent bundle $ T_{\FF(E)/X} $. Therefore, for any such $\mathbf{a}$, we get a corresponding orthogonal decomposition (in the smooth category)
	$$
	T_{\FF(E)} = T_{\FF(E)/X} \oplus T_{\FF(E)/X}^{\perp_{\Xia}}.
	$$
	Finally, from the explicit expression of Formula (\ref{eq: curvature q^a pointwise decreasing}), we see that such a decomposition is independent of the particular choice of the weight $\mathbf{a}$, so that it depends only on $h$ and we can drop any reference to the weight and write
	\begin{equation}\label{eq: splitting tangent projective}
	T_{\FF(E)} = T_{\FF(E)/X} \oplus T_{\FF(E)/X}^{\perp_h}.
	\end{equation}
	We denote by
	\[
	p_1 \colon T_{\FF(E)} \to T_{\FF(E)/X} \ \text{ and } \ p_2 \colon T_{\FF(E)} \to T_{\FF(E)/X}^{\perp_h}
	\]
	the natural projections relative to the splitting \eqref{eq: splitting tangent projective}.
	
	\begin{rem}\label{rem:pvspi}
	Observe that the above splitting is compatible with $d\pi$ in the following sense. The restriction $d\pi|_{T_{\FF(E)/X}^{\perp_h}}\colon T_{\FF(E)/X}^{\perp_h}\to\pi^* T_X$ is a smooth isomorphism of complex vector bundles, and moreover, by a direct pointwise computation in the local holomorphic coordinates chosen above, we have $d\pi=d\pi\circ p_2$, and $d\pi\circ p_1=0$. 
	
	In particular, observe that, for any given point $(x_0,\ff_0)\in\mathbb F_\rho(E)$, given the holomorphic coordinates $(z,\zeta)$ centered at $(x_0,\ff_0)$ as above, we explicitly have that 
$$
\begin{aligned}
& T_{\FF(E)/X,(x_0,\ff_0)}=\Span\{\partial/\partial\zeta_{\lambda\mu}|_{(z,\zeta)=(0,0)}\},\\
&T_{\FF(E)/X,(x_0,\ff_0)}^{\perp_h}=\Span\{\partial/\partial z_k|_{(z,\zeta)=(0,0)}\}.
\end{aligned}
$$
	\end{rem}	
	
The remark above shows that, for $1\le j\le m$, if we define
	\[
	 \Xi_{\rho,j}^{\textrm{vert}}  :=  \Xi_{\rho,j} \circ (p_1\otimes\overline{p_1})
	\]
	and
	\[
	\Xi_{\rho,j}^{\textrm{hor}} :=  \Xi_{\rho,j} \circ (p_2\otimes\overline{p_2}),
	\]
	which are both sections of $ \bigwedge^{1,1}{T_{\FF(E)}^{\vee}} $, we have
	\[
	\Xi_{\rho,j} = \Xi_{\rho,j}^{\textrm{vert}} + \Xi_{\rho,j}^{\textrm{hor}}.
	\]
	Moreover, again with the same choice of coordinates, pointwise at $(x_0,\ff_0)$ we have	
	\begin{equation}\label{eq:hor}
	\bigl(\Xi_{\rho,j}^{\textrm{hor}}\bigr)_{(x_0,\ff_0)} = \frac{i}{2\pi} \sum_{\lambda=s_{j-1} + 1}^{s_j} \Theta_{\lambda\lambda}(x_0).
	\end{equation}
	
	Next, for $ x \in X, $ let $ \ff \in \FF(E_x) $ be given by a unitary basis $ (v_1,\ldots,v_r) $ of $ E_x $.
	Define a section $ \theta_j \colon \FF(E) \to \bigwedge^{1,1}{T_{\FF(E)}^{\vee}} $ as follows
	\begin{equation*}
	\theta_j(x,\ff) = \frac{i}{2\pi} \sum_{\lambda=s_{j-1} + 1}^{s_j} \bigl\langle \pi^*\Theta(E,h)_{(x,\ff)}\cdot {v_\lambda}, {v_\lambda} \bigr\rangle_h.
	\end{equation*}
	\begin{prop}\label{prop:blocks}
		The section $ \theta_j $ is well defined, \textsl{i.e.} it does not depend upon the choice of a particular representative $\mathbf v= (v_1,\ldots,v_r) $ for $\mathbf f$. 
		\begin{proof}
			Take a local normal frame $ (e_1,\ldots,e_r) $ of $ E $ at $ x $ such that $ \mathbf{e} := \bigl(e_1(x),\ldots,e_r(x)\bigr) $ and $ \mathbf{v} $ identify the same flag $ \ff $.
			For $ \lambda = 1,\ldots,r $, we have
			\begin{equation}\label{eq: change of base from e to v}
			v_\lambda = a_{1\lambda} e_1(x) + \cdots + a_{r\lambda}e_r(x)
			\end{equation}
			where $ A = (a_{pq}) $ is the change of coordinates matrix.
			Since the unitary bases $ \mathbf{e} $ and $ \mathbf{v} $ give both the same flag, it follows that $ A $ is a block matrix with the following form
			\[
			\begin{pmatrix}
			A_{11} & 0 & 0 \\
			0 & \ddots & 0 \\
			0 & 0 & A_{mm}
			\end{pmatrix},
			\]
			where the diagonal block $ A_{jj} $ is again a unitary matrix of size $ s_{j} - s_{j-1} $; in particular, for $ 1 \le \lambda \le r $ if the entry $ (\lambda,\lambda) $ hits the block $ A_{jj} $, \textsl{i.e.} if $s_{j-1}<\lambda\le s_j$, Formula \eqref{eq: change of base from e to v} reads as
			\[
			v_\lambda = \sum_{k=s_{j-1}+1}^{s_j}a_{k \lambda} e_{k}(x).
			\]
			Using the local normal frame $\mathbf{e}$ we get
			\begin{align*}
			& \frac{i}{2\pi} \sum_{\lambda=s_{j-1} + 1}^{s_j} {\bigl\langle {{\pi^*\Theta(E,h)}_{(x,\ff)}}\cdot {v_\lambda}, {v_\lambda} \bigr\rangle_h} \\
			&= \frac{i}{2\pi} \sum_{\lambda=s_{j-1} + 1}^{s_j} {\left\langle {{\pi^*\Theta(E,h)}_{(x,\ff)}}\cdot \left( \sum_{\alpha=s_{j-1} + 1}^{s_j} a_{\alpha \lambda} e_\alpha \right), \sum_{\beta=s_{j-1} + 1}^{s_j} a_{\beta \lambda} e_\beta \right\rangle_h} \\
			&= \frac{i}{2\pi} \sum_{\alpha,\beta=s_{j-1} + 1}^{s_j}\underbrace{ \sum_{\lambda=s_{j-1} + 1}^{s_j} a_{\alpha \lambda} \overline{a_{\beta \lambda}}}_{=\delta_{\alpha\beta}} \,\bigl\langle {{\pi^*\Theta(E,h)}_{(x,\ff)}}\cdot {e_\alpha}(x), {e_\beta}(x) \bigr\rangle_h \\
			&= \frac{i}{2\pi} \sum_{\lambda=s_{j-1} + 1}^{s_j}\bigl\langle {{\pi^*\Theta(E,h)}_{(x,\ff)}}\cdot {e_\lambda}(x), {e_\lambda}(x) \bigr\rangle_h,
			\end{align*}
			and the proposition follows.
		\end{proof}
	\end{prop}
	
Of course, $ \theta_j $ is smooth, hence it is a form in $ \mathcal{A}^{1,1}(\FF(E)) $.
	
	\begin{lem}\label{lem: definizione di theta}
		The equality
		\[
		\theta_j = \Xi_{\rho,j}^{\operatorname{hor}}
		\]
		holds.
		\begin{proof}
			At any given $ (x,\ff) \in \FF(E) $, choose $ (e_1,\ldots,e_r) $ to be a local normal frame for $ E $ at $ x $ such that $ \ff $ is given by $\bigl(e_1(x),\ldots,e_r(x)\bigr)$, and consider the induced holomorphic coordinates around $ (x,\ff)$ as above.
			
			Since the evaluation of $ \theta_j $ in $ (x,\ff) $ does not depend on the choice of the unitary basis defining $ \ff $, we have the following chain of equalities
			\begin{align*}
			\theta_j(x,\ff) &= \frac{i}{2\pi} \sum_{\lambda=s_{j-1} + 1}^{s_j} \bigl\langle \pi^*\Theta(E,h)_{(x,\ff)}\cdot {e_\lambda(x)}, {e_\lambda(x)} \bigr\rangle_h\\
			&= \frac{i}{2\pi} \sum_{\lambda=s_{j-1} + 1}^{s_j} \pi^*\Theta_{\lambda\lambda}(x,\ff) = \frac{i}{2\pi} \sum_{\lambda=s_{j-1} + 1}^{s_j} \Theta_{\lambda\lambda}(x)\\
			&= \bigl( \Xi_{\rho,j}^\textrm{hor} \bigr)_{(x,\ff)},
			\end{align*}
			where the last equality follows from Formula \eqref{eq:hor}.
		\end{proof}
	\end{lem}
	
	Now, in order to simplify the notation in what follows, let us relabel $ \Xi_{\rho,j}^{\textrm{vert}}=\omega_j$.
	The section $ \omega_j $ is smooth since, for instance, by Lemma \ref{lem: definizione di theta} it equals $ \Xi_{\rho,j} - \theta_j $, which are two smooth $ (1,1) $-forms on $ \FF(E) $.
	
	Summing up, for $ 1 \le j \le m $, we have shown that
	\begin{equation}\label{eq: first chern form of q_j is sum of two forms}
	\Xi_{\rho,j} = \theta_j + \omega_j, 
	\end{equation}
	that is, the first Chern form of $ (Q_{\rho,j},h) $ can be written as a sum of a horizontal and a vertical part with respect to the natural splitting defined above, where the horizontal part contains the information coming from the curvature of $E$ while the vertical part is the standard curvature of the determinant of the tautological successive quotients on a flag manifold.
	
	\section{Pointwise Gysin's formul\ae{} \textsl{\`{a} la} Darondeau--Pragacz}
	We are now in a good position to prove the following main technical proposition.
	\begin{prop}\label{lem: pol homogeneous}
		Let $ F(\Xi) := F(\Xi_{\rho,1},\ldots,\Xi_{\rho,m}) $ be a complex homogeneous polynomial in the $ (1,1) $-forms $ \Xi_{\rho,1},\ldots,\Xi_{\rho,m} $ on $ \FF(E) $.
		Then the push-forward $ \pis F(\Xi) $ is given by a universal (weighted) homogeneous polynomial evaluated in the Chern forms of $ (E,h) $.
		\end{prop}
		
		In the statement, by universal we mean that this polynomial depends only on the shape of $F$ and on the rank of $E$.
		
		\begin{proof}
			Denote by $ d_{\rho} $ the relative dimension of the flag bundle $ \FF(E) $ and write $ F(\Xi) $ as
			\begin{equation*}
			\sum_{j_1 + \cdots + j_m = d_{\rho}+k} a_{(j_1, \ldots, j_m)} \Xi_{\rho,1}^{j_1} \wedge \ldots \wedge \Xi_{\rho,m}^{j_m},\quad a_{(j_1, \ldots, j_m)}\in\mathbb C,
			\end{equation*}
			where we can w.l.o.g. suppose that $0\le k \le n $, otherwise the push-forward would be identically zero for obvious degree reasons.
			From now on, in order to simplify the notation, we omit the symbol $ \wedge $ for the wedge product of forms and we use, where useful, the multi-index notation.
			Thanks to Formula \eqref{eq: first chern form of q_j is sum of two forms}, we can write
			\begin{align*}
			F(\Xi) &= \sum_{|J| = d_{\rho}+k} a_{J} (\theta_1 + \omega_1)^{j_1} \cdots (\theta_m + \omega_m)^{j_m} \\
			&= \sum_{|J| = d_{\rho}+k} a_{J} \sum_{b_1=0}^{j_1} \binom{j_1}{b_1} \theta_1^{j_1-b_1} \omega_1^{b_1} \cdots \sum_{b_m=0}^{j_m} \binom{j_m}{b_m} \theta_m^{j_m-b_m} \omega_m^{b_m} \\
			&= \sum_{|J| = d_{\rho}+k} a_{J} \sum_{b_1=0}^{j_1} \cdots \sum_{b_m=0}^{j_m} \binom{j_1}{b_1} \cdots \binom{j_m}{b_m} \theta_1^{j_1-b_1} \cdots \theta_m^{j_m-b_m} \omega_1^{b_1} \cdots \omega_m^{b_m}.
			\end{align*}
			
			Since $ \pi $ is a proper submersion, and by definition the push-forward $ \pis F(\Xi) $ is given by integration along the fibers obtained by locally splitting the variables $ (x,\ff) $, at $ x \in X $, we have that
			\[
			{\pis F(\Xi)}_x = \int_{\ff \in \FF(E_x)} {F(\Xi)}_{(x,\ff)},
			\]
			where the second term stands for the integral of $ {F(\Xi)}_{(x,\ff)} $ performed only on those differentials related to the variable $ \ff $. Therefore, for degree reasons, the only terms which can possibly survive after integration along the fibers are those for which $b_1+\cdots+b_m=d_{\rho}$, since the $\omega_j$'s and only the $\omega_j$'s contain the relevant vertical differentials.
			
			For the push-forward we thus obtain
			$$
			\pi_* F(\Xi)_x= \!\!\!\sum_{\substack{|J| = d_{\rho}+k \\ b_1=0,\dots,j_1 \\ \cdots \\ b_m=0,\dots,j_m \\ b_1+\cdots+b_m=d_{\rho}}}\!\!\! a_{J}\binom{j_1}{b_1} \cdots \binom{j_m}{b_m} \int_{\FF(E_x)} \underbrace{\theta_1^{j_1-b_1} \cdots \theta_m^{j_m-b_m}}_{=:\theta^{J-B}}\underbrace{\omega_1^{b_1} \cdots \omega_m^{b_m}}_{=:\omega^B}.
			$$
			
			What we want to do now is to explicitly write $\theta^{J-B}$ at an arbitrary point $(x,\ff)\in\FF(E_x)$. Let $\ff$ be given by a unitary basis $ (v_1,\ldots,v_r)$ of $E_x$, so that $ \theta_\ell(x,\ff) = {i / 2\pi} \sum_{\lambda=s_{\ell-1} + 1}^{s_\ell} {\langle \pi^*\Theta(E,h)_{(x,\ff)}\cdot {v_\lambda}, {v_\lambda} \rangle_h} $. Suppose as above we had fixed a local normal frame $(e_1,\dots,e_r)$ for $E$ centered at $x\in X$, and let $v_{\lambda} = \sum_l v_{\lambda}^l\,e_l(x)$. Thus, we have
			\begin{equation}\label{eq: horizontal form in coordinates}
				\theta_\ell (x,\ff) = \frac{i}{2\pi} \sum_{\alpha,\beta=1}^{r} \left( \sum_{\lambda=s_{\ell-1} + 1}^{s_\ell} v_{\lambda}^{\alpha} {\bar{v}_{\lambda}^{\beta}} \right) \Theta_{\beta \alpha}(x) .
			\end{equation}
	Then, we get the following expression for $\theta^{J-B}$:
	$$
	\begin{aligned}
	\theta^{J-B} &=
	\biggl( \frac{i}{2\pi} \biggr)^{k}
	\bigwedge_{\ell=1}^{m}
	\left(
	\sum_{\alpha_\ell,\beta_\ell = 1}^{r} \left( \sum_{\lambda_\ell = s_{\ell-1} + 1}^{s_\ell} v_{\lambda_\ell}^{\alpha_\ell} {\bar{v}_{\lambda_\ell}^{\beta_\ell}}
	\right)
	\Theta_{\beta_\ell \alpha_\ell}
	\right)^{j_\ell - b_\ell} \\
	&= \biggl( \frac{i}{2\pi} \biggr)^{k}\!\!\!\!\!\!\!\!\!
	\sum_{\substack{
	 		\alpha_1^1, \beta_1^1, \dots, \alpha_1^{j_1-b_1}, \beta_1^{j_1-b_1} = 1,\dots,r \\
	 		\cdots\\
	 		\alpha_m^1, \beta_m^1, \dots, \alpha_m^{j_m-b_m}, \beta_m^{j_m-b_m} = 1,\dots, r}}\!\!\!\!\!\!\!\!\!\!\!\!\!\!\!\!\!\!
		{^BQ}_{\beta_1^1 \cdots \beta_m^{j_m-b_m}}^{\alpha_1^1 \cdots \alpha_m^{j_m-b_m}}
	\Theta_{\beta_1^1 \alpha_1^1} \cdots \Theta_{\beta_m^{j_m-b_m} \alpha_m^{j_m-b_m}},
	\end{aligned}
	$$
	where
	\[
	{^BQ}_{\beta_1^1 \cdots \beta_m^{j_m-b_m}}^{\alpha_1^1 \cdots \alpha_m^{j_m-b_m}} =
	\prod_{\ell=1}^{m}
	\left[
	\left(
	\sum_{\lambda_\ell = s_{\ell-1} + 1}^{s_\ell} v_{\lambda_\ell}^{\alpha_\ell^{1}}\bar{v}_{\lambda_\ell}^{\beta_\ell^{1}}
	\right)
	\cdots
	\left(
	\sum_{\lambda_\ell = s_{\ell-1} + 1}^{s_\ell} v_{\lambda_\ell}^{\alpha_\ell^{j_\ell-b_\ell}}\bar{v}_{\lambda_\ell}^{\beta_\ell^{j_\ell-b_\ell}}
	\right)
	\right].
	\]
	Remark that the $\Theta_{\beta\alpha}$'s only depend on the point $x$, while the $v_\lambda^l$'s can be seen, by a slight abuse of notation, as variables of integration even if they have to be understood modulo the action (cf. with the construction of the matrix $A$ in the proof of Proposition \ref{prop:blocks}) of 
	$$
	U(s_1-s_{0}) \times U(s_{2}-s_{1})\times\cdots \times U(s_m-s_{m-1})\subset U(\underbrace{s_m-s_0}_{=r})
	$$
	which of course corresponds to the homogeneous presentation of the (incomplete) flag manifold as $ U(r) / U(s_1-s_0)\times\cdots \times U(s_m-s_{m-1}) $. 
	
	At the end of the day, the integrals $\int_{\FF(E_x)}\theta^{J-B}\omega^B$ are given by the following sum:
	$$
	 \biggl( \frac{i}{2\pi} \biggr)^{k}
	 \sum_{\substack{
			\alpha_1^1, \ldots, \alpha_m^{j_m-b_m} = 1,\dots, r \\
			\beta_1^1, \ldots, \beta_m^{j_m-b_m} = 1,\dots, r }}
	 {^Bq}_{\beta_1^1 \cdots \beta_m^{j_m-b_m}}^{\alpha_1^1 \cdots \alpha_m^{j_m-b_m}}
	 \Theta_{\beta_1^1 \alpha_1^1} \cdots \Theta_{\beta_m^{j_m-b_m} \alpha_m^{j_m-b_m}},
	$$
	where
	$$
	{^Bq}_{\beta_1^1 \cdots \beta_m^{j_m-b_m}}^{\alpha_1^1 \cdots \alpha_m^{j_m-b_m}} =
	\int_{\FF(E_x)}
	{^BQ}_{\beta_1^1 \cdots \beta_m^{j_m-b_m}}^{\alpha_1^1 \cdots \alpha_m^{j_m-b_m}}
	\,\omega^B.
	$$
	The good news is that these coefficients ${^Bq}_{\beta_1^1 \cdots \beta_m^{j_m-b_m}}^{\alpha_1^1 \cdots \alpha_m^{j_m-b_m}}$ are universal, in the following sense: they do not depend anymore on the metric $h$, nor on the point $x\in X$, but only on the multi-index $B$ and on the rank $r$ of $E$. 
	
	Indeed, they might be calculated in the \lq\lq absolute\rq\rq{} case of the flag manifold $\FF(\mathbb C^r)$, where $\mathbb C^r$ is endowed with the standard Euclidean metric, the $v_\lambda^l$'s are the element of a matrix in $U(r)$ representing the given flag, and where the top form $\omega^B$ is nothing else than the corresponding wedge product of the curvature forms of the determinant of the tautological quotients line bundles on $\FF(\mathbb C^r)$ with respect to the natural metrics induced by the Euclidean metric of $\mathbb C^r$.
			
			Summing up, we have shown that
			$$
			\begin{aligned}
			{\pis F(\Xi)}_x &=
			 \sum_{\substack{|J| = d_{\rho}+k \\ b_1=0,\dots,j_1 \\ \cdots \\ b_m=0,\dots,j_m \\ b_1+\cdots+b_m=d_{\rho}}}
			\sum_{\substack{
			\alpha_1^1, \ldots, \alpha_m^{j_m-b_m} = 1,\dots, r \\
			\beta_1^1, \ldots, \beta_m^{j_m-b_m} = 1,\dots, r }}
			a_{J}\binom{j_1}{b_1} \cdots \binom{j_m}{b_m}{}^Bq_{\beta_1^1 \cdots \beta_m^{j_m-b_m}}^{\alpha_1^1 \cdots \alpha_m^{j_m-b_m}} \\
			&\qquad\qquad\times\left( \frac{i}{2\pi} \Theta_{\beta_1^1 \alpha_1^1} \right) \cdots \left( \frac{i}{2\pi} \Theta_{\beta_m^{j_m-b_m} \alpha_m^{j_m-b_m}} \right) .
			\end{aligned}
			$$
The above expression is thus given by evaluating a homogeneous polynomial $ \tilde P $ of degree 
			$$
			 (j_1-b_1)+(j_2-b_2)+\cdots+(j_m-b_m) =|J|-|B|=k
			 $$ 
on the entries $\Theta_{\beta\alpha}$'s of the matrix associated to the curvature $\Theta(E,h)_x$ with respect to the frame $\bigl(e_1(x),\dots,e_r(x)\bigr)$. 

This polynomial (\textsl{i.e.} its coefficients) is clearly independent of the point $x\in X$. Moreover, $\tilde P$ is of course also invariant under change of frame at $x$, because from its very definition in terms of push-forwards it is independent of the local frame chosen to make the computation. It follows, say from Chern--Weil theory (see for instance \cite[p. 402]{GH78}), that there exists a weighted homogeneous polynomial $P$ such that globally on $X$
\begin{equation}\label{eq: polynomial chern forms pi star xi}
		\pis F(\Xi) = P\bigl(c_1(E,h),\ldots,c_r(E,h)\bigr),
		\end{equation}
		where
		$$ c_j(E,h)  = \Tr_{\Lambda^j E} \left( {\bigwedge}^j \frac{i}{2\pi} \Theta(E,h) \right) $$
		is the $ j $-th \emph{Chern form} of $ (E,h) $. 
By construction, $P$ is universal since it obviously depends only on $F$ and the rank of $E$.	
		\end{proof}
	
	\begin{rem}\label{rem:chernvssegre}
		Of course, since Chern forms may be expressed in terms of Segre forms, a completely analogous statement holds with a polynomial whose variables are now the Segre forms $s_j(E,h)$ of $(E,h)$.
	\end{rem}
	
	\subsection{Main result}
	We shall see now how to compute the push-forward of the differential form $ F(\Xi) $ through the projection $ \pi $, where $ F $ is a complex homogeneous polynomial in $ \Xi_{\rho,1}, \ldots, \Xi_{\rho,m} $.
	In order to do this, we use a formula by Darondeau and Pragacz given in \cite[Proposition 1.2]{DP17}, which allows to compute the push-forward at the level of cohomology (in fact at the level of Chow rings, but cohomology suffices for our purposes). 
	
	First of all, observe that, even if it is stated for algebraic manifolds over an algebraically closed field $k$, when $k=\mathbb C$ it also holds for non necessary algebraic complex manifolds. Indeed, ultimately, their proof only relay upon the push-forward formula for the tautological class on the projectivized bundle in terms of Segre classes, which is definitely valid also for general complex manifolds.
	
	So, from now on, we suppose that $ E $ is a rank $ r $ holomorphic vector bundle over a complex $ n $-dimensional manifold $ X $ and denote by $ d_{\rho} $ the relative dimension of $ \FF(E) $.
	In what follows, we use the same symbol for the push-forward
	$$ \pis \colon \mathcal{A}^{2(d_{\rho}+k)}(\FF(E)) \to \mathcal{A}^{2k}(X) $$
	of differential forms and the one induced in cohomology, namely
	$$ \pis \colon H^{2(d_{\rho}+k)}(\FF(E)) \to H^{2k}(X) . $$
	With the same notation of \cite{DP17}, we denote by $ \xi_1, \ldots, \xi_r $ the (virtual) Chern roots of $ \pi^{*} E^{\vee} $.
	The Darondeau--Pragacz formula allows us to compute the push-forward of any cohomology class of the form
	$$ \tilde F(\xi_1,\ldots,\xi_r) \in H^{2(d_{\rho}+k)}(\FF(E)) $$
	in terms of the Segre classes $ s_j := s_j(E), 1\le j \le n $, of the vector bundle $ E $ (we are thus implicitly asking here that the polynomial $\tilde F$ has the appropriate symmetries).
	
	More precisely, let $ t_1,\ldots,t_r $ be a set of formal variables, then
	\begin{equation}\label{eq: darondeau-pragacz classes coefficient}
	\begin{aligned}
	&\pis \tilde F(\xi_1,\ldots,\xi_r) \\
	&= [t_1^{\ell_1} \cdots t_r^{\ell_r}] \left( \tilde F(t_1,\ldots,t_r) \prod_{1\le i \le r} \left( 1 + \sum_{j=1}^{n} \frac{s_j}{t_i^j} \right) \prod_{1\le i < j \le r} (t_i - t_j) \right)
	\end{aligned}
	\end{equation}
	where, for a monomial $ \mathfrak m $ and a (Laurent) polynomial $ P $, the notation $ [\mathfrak m](P) $ stands for the coefficient of $ \mathfrak m $ in $ P $, and the rule to determine the $\ell_j$'s is as follows: for $r-\rho_k<j\le r-\rho_{k-1}$, say $j=r-\rho_k+i$ for some $i=1,\dots,\rho_k-\rho_{k-1}$, we set $\ell_j=r-i$.
This is (one possible instance of) Darondeau--Pragacz formula.

	\begin{rem}\label{rem: universality of phi}
		Let us call $ \Phi(s_1,\ldots,s_n) $ the right-hand side of Formula \eqref{eq: darondeau-pragacz classes coefficient}.
		If we consider the Segre classes $ s_1,\ldots,s_n $ as formal variables, we can affirm that the polynomial $ \Phi(s_1,\ldots,s_n) $ is universal, in the sense that its coefficients depends only upon $ \tilde{F} $ and the rank $ r $.
		Moreover, $ \Phi $ is, by construction, weighted homogeneous of degree $ 2k $, since $ \deg s_j = 2j $.
	\end{rem}

\begin{rem}\label{rem:symmetries}
Of course, for our purposes, it shall suffice to use some of the possible polynomials $\tilde F$ considered above only, namely those whose symmetries are of the form
$$
\tilde F(\xi_1,\dots,\xi_r)=F\left(\dots,-\sum_{s_{j-1}+1}^{s_j}\xi_\ell,\dots\right)=F\bigl(c_1(Q_{\rho,1}),\dots,c_1(Q_{\rho,m})\bigr).
$$
Observe however, that in the special case of \emph{complete} flag bundles, the classes $\xi_1,\dots,\xi_r$ are not virtual and give actual cohomology classes, so that no further requirements for symmetries of $\tilde F$ are needed.
\end{rem}	
	
	The key result of this section is that, in fact, Darondeau--Pragacz formula~\eqref{eq: darondeau-pragacz classes coefficient} also holds pointwise at the level of differential forms, in the Hermitian setting, for polynomials $\tilde F$ satisfying the further symmetries considered in Remark \ref{rem:symmetries}.
	
	\begin{thm}\label{thm: darondeau-pragacz for differential forms}
		Let $ (E,h) $ be a rank $ r $ Hermitian holomorphic vector bundle over a complex manifold $ X $ of dimension $ n $, and let $ F $ be a complex homogeneous polynomial of degree $ d_{\rho}+k $ in $m$ variables.
		Then, we have the equality
		\begin{align*}
		\pis F(\Xi_{\rho,1},\ldots,\Xi_{\rho,m}) = \Phi\bigl(s_1(E,h),\ldots,s_n(E,h)\bigr).
		\end{align*}
		\end{thm}
Clearly, in the statement above, $ \Phi $ is the polynomial introduced in Remark \ref{rem: universality of phi} associated to the polynomial $\tilde F(\xi_\bullet)=F\bigl(c_1(Q_{\rho,\bullet})\bigr)$, as in Remark \ref{rem:symmetries}. 
		\begin{proof}
			If $ \eta \in \mathcal{A}^{\bullet}(X) $ is a closed form, then $ [\eta] $ stands as usual for the cohomology class in $ H^{\bullet}(X) $ represented by $ \eta $.
			By Formula \eqref{eq: darondeau-pragacz classes coefficient}, it holds that
			\begin{align*}
			[\pis F(\Xi_{\rho,1},\ldots,\Xi_{\rho,m})] &= \pis \bigl[ F(\Xi_{\rho,1},\ldots,\Xi_{\rho,m})\bigr] \\
			&= \pis \tilde F(\xi_1,\ldots,\xi_r) \\
			&=  \Phi(s_1,\ldots,s_n) \\
			&= \left[ \Phi\bigl(s_1(E,h),\ldots,s_n(E,h)\bigr) \right].
			\end{align*}
			Hence, the difference
			\begin{align*}
			\pis F(\Xi_{\rho,1},\ldots,\Xi_{\rho,m}) - \Phi(s_1(E,h),\ldots,s_n(E,h))
			\end{align*}
			must be an exact global $ (k,k) $-form on $ X $.
			Recall that by Proposition \ref{lem: pol homogeneous} and Formula \eqref{eq: polynomial chern forms pi star xi}, $\pis F(\Xi_{\rho,1},\ldots,\Xi_{\rho,m}) $ is a universal weighted homogeneous polynomial $ P=P\bigl(c_\bullet(E,h)\bigr) $ of weighted degree $ 2k $ in the Chern forms of $ (E,h) $.
			If we express the Segre forms in terms of the Chern forms, the previous difference can be written as a complex weighted homogeneous polynomial
			\begin{equation*}
			G\bigl(c_1(E,h),\ldots,c_r(E,h)\bigr) = \sum_{k_1 + 2k_2 +\cdots + rk_r = k} g_{k_1 \cdots k_r} c_1(E,h)^{k_1} \wedge \cdots \wedge {c_r(E,h)}^{k_r}
			\end{equation*}
			in the Chern forms.
			Note that $ G $ is universal (since $ P $ and $ \Phi $ are) in the sense that its coefficients $ g_{k_1 \cdots k_r} $ do not depend upon $ E $, nor $X$, but only upon $r,n$, and $F$. Recall that our aim is to show that $G$ is in fact identically zero: following \cite{Gul12}, to achieve this we shall evaluate it on a particular vector bundle on a particular class of  manifolds, as follows.
			
			Take $X$ to be any $n$-dimensional \emph{projective} manifold and fix an ample line bundle $ A $ on $ X $. Let $ \omega_A $ be a metric on $ A $ with positive curvature.
			For $ m_1,\ldots, m_r $ positive integers, we define the totally split, rank $r$ vector bundle
			\[
			\mathcal{E} := {A}^{\otimes m_1} \oplus \cdots \oplus {A}^{\otimes m_r},
			\]
			and denote by $ \omega_{\mathcal{E}} $ the natural induced metric on $ \mathcal{E} $ by $ \omega_A $.
			Hence,
			\begin{align*}
			G&(c_1(\mathcal{E},\omega_{\mathcal{E}}),\ldots,c_r(\mathcal{E},\omega_{\mathcal{E}})) \\
			&= \sum_{k_1 + 2k_2 +\cdots + rk_r = k} g_{k_1 \cdots k_r} c_1(\mathcal{E},\omega_{\mathcal{E}})^{k_1} \wedge \cdots \wedge {c_r(\mathcal{E},\omega_{\mathcal{E}})}^{k_r} \\
			&= \sum_{k_1 + 2k_2 +\cdots + rk_r = k} g_{k_1 \cdots k_r} \bigwedge_{s=1}^{r} \left( {\sum_{1 \le j_1 < \cdots < j_s\le r}} m_{j_1}\cdots m_{j_s} {c_1(A,\omega_A)}^{s} \right)^{k_s} \\
			&= \left[ \sum_{k_1 + 2k_2 +\cdots + rk_r = k} g_{k_1 \cdots k_r} \prod_{s=1}^{r} \left( {\sum_{1 \le j_1 < \cdots < j_s\le r}} m_{j_1}\cdots m_{j_s} \right)^{k_s} \right] {c_1(A,\omega_A)}^{k} .
			\end{align*}
			Let $ T_1, \ldots, T_r $ be a set of formal variables, and consider the polynomial $ p $ defined by
			\begin{equation}\label{eq: forma di p}
			p(T_1,\ldots,T_r) = \sum_{k_1 + 2k_2 +\cdots + rk_r = k} g_{k_1 \cdots k_r} \prod_{s=1}^{r} \left( {\sum_{1 \le j_1 < \cdots < j_s\le r}} T_{j_1}\cdots T_{j_s} \right)^{k_s} .
			\end{equation}
			We have by definition that
			\[
			G(c_1(\mathcal{E},\omega_{\mathcal{E}}),\ldots,c_r(\mathcal{E},\omega_{\mathcal{E}})) = p(m_1,\ldots,m_r) {c_1(A,\omega_A)}^{k}
			\]
			and, consequently, in cohomology it holds that
			\[
			p(m_1,\ldots,m_r) c_1(A)^k = \left[ G(c_1(\mathcal{E},\omega_{\mathcal{E}}),\ldots,c_r(\mathcal{E},\omega_{\mathcal{E}})) \right] = 0 .
			\]
			Since $ A $ is ample, the only possibility is that the polynomial $ p $ is zero for every choice of positive integers $ m_1,\ldots,m_r $. But the set of points in $\mathbb C^r$ whose coordinates are positive integers is Zariski dense, and thus $p$ must be identically zero. Consequently, having the same coefficients as $p$, $ G \equiv 0 $ and this concludes the proof.
\end{proof}

			
			
	\begin{rem}\label{rem:schurfunctions}
		As already said, Formula \eqref{eq: darondeau-pragacz classes coefficient} is not the only possible instance of the Gysin formula given by Darondeau--Pragacz. For example, they gives in \cite[Proposition 4.2]{DP17} the universal Gysin formula for flag bundles in terms of \emph{Schur functions} (a sort of generalization of Schur polynomials, see next section for a definition of Schur polynomial).
		
		Such version of the Darondeau--Pragacz formula is particularly useful since for instance it explicitly shows (with a little further manipulation, see the forthcoming second-named author PhD thesis, as well as \cite{Fag21}, for more details) that one can obtain \emph{all} the Schur polynomials in the Chern classes of $E$ as a push-forward from the complete flag bundle of monomials of type $ (-\xi_1)^{\lambda_1} \cdots (-\xi_r)^{\lambda_r} $, provided the $ \lambda_j$'s,\ satisfy a certain relation (for more details see again \cite[Proposition 4.2]{DP17}).
		
		Clearly, the validity of \cite[Proposition 4.2]{DP17} at the level of differential forms follows directly from our Theorem \ref{thm: darondeau-pragacz for differential forms}.
	\end{rem}

	\section{Application to Griffiths' conjecture}\label{sec:positivity}
	Let us first recall a few notations about Schur polynomials, essentially taken from the exposition in \cite[\S 2]{DPS94} (see also \cite{FL83}).
	
	Denote by $ \Lambda(k,r) $ the set of all the partitions $ \sigma = (\sigma_1,\ldots,\sigma_k) $ in $ \mathbb N^{k} $ such that
		\[
		r \ge \sigma_1 \ge \ldots \ge \sigma_k \ge 0, \quad |\sigma|=\sum_{j=1}^{k} \sigma_j = k .
		\]
		For every $ \sigma \in \Lambda(k,r) $ we can define a \emph{Schur polynomial} $ S_{\sigma} \in \Z[c_{1},\ldots,c_{r}] $ of weighted degree $ 2k $ (where $ \deg c_{j} = 2j $) as
		\[
		S_{\sigma}(c_{1},\ldots,c_{r}) := \det
		\begin{pmatrix}
		c_{\sigma_1} & c_{\sigma_1+1} & \cdots & c_{\sigma_1+k-1} \\
		c_{\sigma_2-1} & c_{\sigma_2} & \cdots & c_{\sigma_2+k-2} \\
		\vdots & \vdots & \ddots & \vdots \\
		c_{\sigma_k-k+1} & c_{\sigma_k-k+2} & \cdots & c_{\sigma_k}
		\end{pmatrix}
		\]
		where, by convention, $ c_0 = 1 $ and $ c_j = 0 $ if $ j \notin [0,r] $.
		
		The Schur polynomials, as $\sigma\in\Lambda(k,r)$ varies, form a basis for the $ \Q $-vector space of degree $ 2k $ weighted homogeneous polynomials in $ r $ variables.
		Thus, given such a polynomial $ P $ we can write
		\begin{equation*}
		P = \sum_{\sigma \in \Lambda(k,r)} b_{\sigma}(P) S_{\sigma}.
		\end{equation*}
		The set of all $P$ such that $b_{\sigma}(P) \ge 0$ for every $\sigma \in \Lambda(k,r)$, which is called the set of \emph{positive polynomials}, is of course a positive convex cone, which we call $\Pi(r)$ following \cite{Gri69} (remark that this is not exactly the positive cone considered by Griffiths, but they coincides \textsl{a posteriori} thanks to the work of \cite{FL83}). It is well known that any product of Schur polynomials can be written as a linear combination of Schur polynomials with non-negative integral coefficients; the values of these coefficients is given combinatorially by the Littlewood--Richardson rule. Thus, these positive cones are stable under product (cf. with the analogous property for wedge product of strongly positive forms, as observed in Remark \ref{rem:stronglyGriffiths}).
		
		Now, if $ E $ is a rank $r$ holomorphic vector bundle over the complex manifold $ X $ and if $ \sigma $ is a partition in $ \Lambda(k,r) $, the \emph{Schur class} of $ E $ associated to $ \sigma $ is the cohomology class
		\begin{equation*}
S_{\sigma}(E) := S_{\sigma}\bigl(c_{1}(E),\ldots,c_{r}(E)\bigr)\in H^{2k}(X,\mathbb Z)
\end{equation*}
formally obtained by computing $S_\sigma$ on the the Chern classes of $E$.

		In the same way, if $ (E,h) $ is a holomorphic Hermitian rank $r$ vector bundle on $ X $, we can define the \emph{Schur form} of $ (E,h) $ associated to $ \sigma $ formally obtained by computing $S_\sigma$ on the the Chern forms of $(E,h)$, and we denote it by $ S_{\sigma}(E,h) $.
		Clearly, the closed differential $2k$-form $S_{\sigma}(E,h)$ is a special representative for the class $ S_{\sigma}(E) $.
		\begin{ex}
		By definition,  the $ k $-th Chern class of $ E $ corresponds to the partition $ (k,\underbrace{0,\dots,0}_{\textrm{$k-1$ times}}) \in\Lambda(k,r)$, \textsl{i.e.}
		\[
		S_{(k,0,\ldots,0)}(E) = c_k(E),
		\]
		while the partition $ (\underbrace{1,\dots,1}_{\textrm{$k$ times}})\in\Lambda(k,r)$ give rise to the $k$-th signed Segre class
		\begin{equation*}
		S_{(1,\ldots,1)}(E) = (-1)^k s_k(E).
		\end{equation*}
		Clearly, at the level of differential forms we get similar equalities.
	\end{ex}
	
As we saw in the introduction, Griffiths conjectured (and proved partially) in \cite{Gri69} that given any rank $r$ Hermitian holomorphic positive vector bundle on a projective manifold, the polynomials belonging to $\Pi(r)$ whenever evaluated on its Chern classes have to return a positive number once integrated over any subvariety of the correct dimension. A full proof of this conjecture is given in \cite{FL83} in the more general setting of ample vector bundles (see also \cite{DPS94} for the even more general context of $E$ nef and $X$ compact K\"ahler).

Actually, in \cite{FL83} a more universal problem is considered and settled, \textsl{i.e.} to characterize precisely  the \emph{numerically positive polynomials for ample vector bundles of rank $r$}. These are defined to be those weighted homogeneous polynomials say of degree $2n$ that whenever evaluated on the Chern classes of any rank $r$ ample vector bundle $E$ over an irreducible projective variety of dimension $n$ give a positive number. Once again, the characterization is that numerically positive polynomials for ample vector bundles of rank $r$ are exactly the non zero positive polynomials.

\begin{rem}\label{rem:smooth vs singular}
It is observed in \cite{FL83}, in Remark (1) right after the proof of Proposition 3.4, that if a weighted homogeneous polynomial $P$ of degree $2n$ is not positive, \textsl{i.e.} it does not belong to $\Pi(r)$, then there exists a \emph{smooth} projective manifold of dimension $n$, and an ample vector bundle of rank $r$ over it, such that when one evaluates this polynomial in its Chern classes and integrates over the manifold, one gets a negative number. Moreover, such vector bundle is constructed as a quotient of a direct sum of very ample line bundles.

The upshot is that if we want to show that a weighted homogeneous polynomial is positive it suffices to show that it is a numerically positive polynomial for ample vector bundles over \emph{smooth} projective manifolds. This will be useful later during the proof of Theorem \ref{thm: positivity of push-forwards by flag bundles}.
\end{rem}

Now, recall that given a Hermitian holomorphic vector bundle $(E,h)\to X$, it is said to be \emph{Griffiths semipositive} (resp. \emph{Griffiths positive}) if for every $x\in X$, $v\in E_x$, $\tau\in T_{X,x}$ we have
$$
\langle\Theta(E,h)_x\cdot v,v\rangle_h(\tau,\bar\tau)\ge 0
$$
(resp. $ > 0 $ and $=0$ if and only if $v$ or $\tau$ is the zero vector). If we compute the Chern curvature at the given point $x$ with respect to a unitary frame for $E_x$, say 
$$
\Theta(E,h)_x= \sum_{\alpha,\beta=1}^r \sum_{p,q=1}^n c_{pq\alpha\beta}(x)\,dz_p \wedge d\bar{z}_q \otimes e_{\alpha}^{\vee} \otimes e_{\beta},
$$
then this is equivalent to ask the same inequalities for the quantity
$$
\sum_{\alpha,\beta} \sum_{p,q} c_{pq\alpha\beta}(x)\,\tau_p \bar\tau_q v_\alpha\bar v_\beta,
$$
where the $v_\alpha$'s and $\tau_p$'s are the coordinate of the chosen vector with respect to the given bases.

A Griffiths positive vector bundle on a compact complex manifold is ample (the converse is not known in general, but it is a conjecture), and a globally generated vector bundle can be endowed with a Hermitian metric which makes it Griffiths semipositive.

It is then natural to ask whether in this Hermitian setting the Fulton--Lazarsfeld--Demailly--Peternell--Schneider theorem holds pointwise for Chern forms, and this is also a question raised by Griffiths in the same paper. Recall that a $(k,k)$-form $u$ is positive if and only if its restriction to every $k$-dimensional complex submanifold is a non negative volume form (cf. \cite[Chapter III, \S 1.A, (1.6) Criterion]{Dem01}). 

\begin{quest}[{\cite{Gri69}}]\label{quest:Griffiths}
Given a Griffiths (semi)positive Hermitian holomorphic vector bundle $(E,h)$, is it true that the positive polynomials evaluated on the Chern forms of $E$ give rise to positive forms?
\end{quest}

\begin{rem}
Coming back to Remark \ref{rem:smooth vs singular}, we see that given $P$ a weighted homogeneous polynomial of degree $2n$ which is not positive, there exists a rank $r$ holomorphic Hermitian vector bundle $(E,h)$ over a smooth projective manifold $X$ of dimension $n$ whose Chern curvature is Griffiths (as well as dual Nakano) positive and such that the corresponding characteristic form obtained by computing $P$ in the Chern forms of $(E,h)$ is not a positive (volume) form. 

This is because one can endow $E$ with the quotient metric of a positively curved direct sum metric on the direct sum of the very ample line bundles in question. Such a metric, begin a quotient of a positively curved (in any sense) metric, is both Griffiths and dual Nakano positive (but not Nakano positive, in general). The corresponding volume form $P\bigl(c_\bullet(E,h)\bigr)$ has negative total mass, and hence must be negative somewhere. 

This means that, even in the pointwise Hermitianized case considered by Griffiths, the cone of positive polynomials is the largest possible for which one can hope such a result.
\end{rem}

Griffiths, in \textsl{loc. cit.}, answered in the affirmative to this question in the special case of the second Chern form of a rank $2$ Griffiths positive holomorphic vector bundle (for the first Chern form the answer is trivially yes). Remark that this question gives, under the stronger hypothesis of Griffiths positivity, a stronger answer than its cohomological version stated earlier, since ---as observed--- a positive polynomial in the Chern form is a special representative in cohomology of the corresponding positive polynomial in the Chern classes. 
	
In the last recent years there have been several partial results towards a fully affirmative answer to Griffiths' question. First, \cite{Gul12} (see also \cite{Div16} for a more direct proof of the main technical result needed, as well as \cite{Mou04} for similar, and somehow more general, computations) proved that the answer is affirmative in the special case of singed Segre classes. 

Then, \cite{Li20} proved the full statement but under the stronger assumption of Bott--Chern (semi)positivity for $(E,h)$. We refer to \textsl{ibid.} for the definition of this variant of Hermitian positivity, which has been observed to be indeed equivalent to dual Nakano (semi)positivity by Finski \cite{Fin20}. 

Other related interesting results about finding positive representatives (not necessarily coming from the given positively curved metric) of the Schur polynomials in the Chern classes are obtained in \cite{Pin18,Xia20}. 

For an even more recent result in the case of (dual) Nakano positive vector bundles, see \cite{Fin20} and Remark \ref{rem:finski}.

\medskip

Here, we are concerned with the original Question \ref{quest:Griffiths}, which is still very much open: one can construct indeed (local, say over a ball) examples of Hermitian holomorphic vector bundles which are Griffiths positive but not Nakano nor dual Nakano positive, see for instance \cite[Proposition 2.9]{Fin20}. Our Theorem \ref{thm: darondeau-pragacz for differential forms} allows us indeed to confirm the strong positivity of quite a few new positive combinations of Schur polynomials in the Chern forms, as follows.

Let $(E,h)$ be a Griffiths semipositive vector bundle over the complex manifold $X$. Consider the flag bundle $ \pi \colon \FF(E) \to X $ and, for $ \mathbf{a} \in \mathbb N^r $ satisfying Condition \eqref{eq:condition on a}, let $ \Qa \to \FF(E) $ be the natural line bundle introduced before. The first observation is contained in the following.

\begin{prop}[{\cite[Lemma 3.7 (a), and Formula (4.9)]{Dem88}}]\label{prop: q^a is positive}
		If $ \mathbf{a} = (a_1,\ldots,a_r)\in\mathbb N^r $ is non increasing, then $ \Qa \to \FF(E) $ endowed with the natural induced Hermitian metric is a semipositive line bundle.
\end{prop}

This means precisely that the Chern curvature $\Xia$ is a closed positive $(1,1)$-form. 

Let us recall that a $(k,k)$-form is \emph{strongly positive} if and only if all of its wedge products against a positive form of complementary bi-degree give a non negative volume form (cf. \cite[Chapter III, (1.1) Definition]{Dem01} for the definition of strongly positive form, the one just given being usually a characterization). Now, strongly positive forms are positive but the converse is not true in general. However, strongly positive $(k,k)$-forms and positive $(k,k)$-forms do coincide for $k=0,1,n-1,n$, \cite[Chapter III, (1.9) Corollary]{Dem01}. Thus, $\Xia$ is also a strongly positive $(1,1)$-form.

Since the wedge product of strongly positive forms is again strongly positive \cite[Chapter III, (1.11) Proposition]{Dem01}, then all wedge powers of $\Xia$ are again strongly positive. Now, it is straightforward to see that the push-forward of a closed strongly positive form under a proper holomorphic submersion is again a closed strongly positive form \cite[Chapter III, (1.17) Proposition]{Dem01}. Thus, we obtain immediately the next proposition.

\begin{prop}\label{prop:pospush}
If $(E,h)\to X$ is a Griffiths semipositive vector bundle, then the closed forms
$$
\pi_*(\Xia)^{d_{\rho}+k},\quad a_1 \ge a_2 \ge \cdots \ge a_r \ge 0,
$$
where $d_{\rho}$ is the relative dimension and $k$ is a non negative integer, are closed strongly positive $(k,k)$-forms.
\end{prop}

\begin{rem}
If it happens that the chain of inequalities in the above statement is not strictly decreasing where prescribed by Condition \eqref{eq:star}, then the push-forward is identically zero. This is because in this case the curvature $\Xia$ has some vertical zero eigenvalue in each fiber, thanks to Formula (\ref{eq: curvature q^a pointwise decreasing}). Therefore, the vertical top form against which we integrate to obtain the push-forward is identically zero being, modulo a factor, the determinant of the vertical part of the curvature. So, in what follows, we can consider without loss of generality only weights $\mathbf a\in\mathbb N^r$ satisfying Condition \eqref{eq:condition on a} and such that $a_{s_1}>a_{s_2}>\cdots >a_{s_m}$.
\end{rem}
	
Now, we come to the main result of this section. Denote by $ \Phi_{\mathbf{a}}^k(E,h) $ the push-forward
$$
\pis(\Xia)^{d_{\rho}+k} = \Phi_{\mathbf{a}}^k(E,h).
$$
It is a closed $(k,k)$-form representing the cohomology class $[\pis(\Xia)^{d_{\rho}+k}]\in H^{2k}(X)$.
	 
	\begin{thm}\label{thm: positivity of push-forwards by flag bundles}
		Let $(E,h)\to X$ be a rank $r$ Griffiths semipositive vector bundle. For every $ \mathbf{a} \in \mathbb N^r $ satisfying Conditions \eqref{eq:condition on a} and \eqref{eq:star}, and for every $ k = 0, \dots, n=\dim X $, the differential form $ \Phi_{\mathbf{a}}^k(E,h) $ is a closed strongly positive $ (k,k) $-form on $ X $ belonging to the positive convex cone $\Pi(r)$ spanned by the Schur forms of $ (E,h) $.

Moreover, the explicit expression of $\Phi_{\mathbf{a}}^k(E,h)$ can be obtained by formally evaluating in the Segre forms $s_j(E,h)$'s of $(E,h)$ the right hand side of Formula (\ref{eq: darondeau-pragacz classes coefficient}), with
$$
\tilde F(t_1,\dots,t_r)=\left(-\sum_{j=1}^m\sum_{\lambda=s_{j-1}+1}^{s_j}a_{s_j}t_\lambda\right)^{d_{\rho}+k}.
$$
\end{thm}	

This theorem covers in particular Guler's work \cite{Gul12}, which concerned push-forwards from the projectivized bundle. 

\begin{rem}\label{rem:stronglyGriffiths}
Observe that it is in some sense more natural to obtain that these forms are strongly positive rather than merely positive. This is because, as said earlier, positive polynomials are stable under products and so do strongly positive forms, while a product of positive forms is not necessarily still positive.
\end{rem}

		\begin{proof}
By Remark \ref{rem:chernvssegre}, there exists a unique homogeneous polynomial $\Gamma^k_{\mathbf a}$ of weighted degree $2k$ such that
$$
\Phi_{\mathbf{a}}^k(E,h)=\Gamma^k_{\mathbf a}\bigl(s_1(E,h),\dots,s_n(E,h)\bigr).
$$
			
			By definition, we have that
			\[
			(\Xia)^{d_{\rho}+k} = \left( a_{s_1} \Xi_{\rho,1} + \cdots + a_{s_m} \Xi_{\rho,m} \right)^{d_{\rho}+k},
			\]
			and thanks to Theorem \ref{thm: darondeau-pragacz for differential forms} we get the explicit expression claimed at the end of the statement. The closedness and strong positivity of $\Phi_{\mathbf{a}}^k(E,h)$ are the content of Proposition \ref{prop:pospush}.

We now want to show that $\Phi_{\mathbf{a}}^k(E,h)$ can be written as a positive linear combination of Schur forms. To do this, let $\tilde\Gamma_{\mathbf{a}}^k$ be the \emph{unique} polynomial such that 
$$
\Phi_{\mathbf{a}}^k(E,h)=\Gamma_{\mathbf{a}}^k\bigl(s_1(E,h),\dots,s_n(E,h)\bigr)=\tilde\Gamma_{\mathbf{a}}^k\bigl(c_1(E,h),\dots,c_r(E,h)\bigr).
$$
Observe that obviously the polynomial $\tilde\Gamma_{\mathbf{a}}^k$ does not depend on the particular vector bundle considered, nor on the particular given base manifold, as usual.
What we want to show will then follow from the Fulton--Lazarsfeld theorem \cite{FL83} if we can prove that $\tilde\Gamma_{\mathbf{a}}^k$ is a numerically positive polynomial for ample vector bundles of rank $r$ over smooth projective manifolds, see Remark \ref{rem:smooth vs singular}.

So, take any rank $r$ ample vector bundle $\mathcal V$ over a $k$-dimensional projective manifold $Z$. By \cite[Lemma 4.1]{Dem88amp}, the corresponding line bundle $\Qa$ over $\FF(\mathcal V)$ is ample, and it can be therefore endowed with a smooth Hermitian metric $h_{\mathcal V,\mathbf a}$ whose Chern curvature $i\Theta(\Qa,h_{\mathcal V,\mathbf a})$ is \emph{strictly} positive, \textsl{i.e.} a K\"ahler form. But then, $\pi_*c_1(\Qa,h_{\mathcal V,\mathbf a})^{d_{\rho}+k}$, is a closed  positive \emph{nowhere zero} $(k,k)$-form representing the cohomology class $\tilde\Gamma_{\mathbf{a}}^k\bigl(c_1(\mathcal V),\dots,c_r(\mathcal V)\bigr)$. In particular, being represented by a non zero positive $(k,k)$-form, we have that
$$
\int_Z\tilde\Gamma_{\mathbf{a}}^k\bigl(c_1(\mathcal V),\dots,c_r(\mathcal V)\bigr)> 0,
$$
as desired.
		\end{proof}


As a byproduct of the proof above one immediately obtains the following statement for ample vector bundles in the same spirit of \cite{Pin18,Xia20}.

\begin{thm}
Let $E\to X$ be an ample vector bundle of rank $r$ over a projective manifold. For every $ \mathbf{a} \in \mathbb N^r $ satisfying Conditions \eqref{eq:condition on a} and \eqref{eq:star}, and for every $k = 0, \dots, n=\dim X$, the $(k,k)$-cohomology classes $\pi_*c_1(\Qa)^{d_{\rho}+k}$ contain a closed strongly positive form and belong to $\Pi(r)$ .
\end{thm}

\begin{rem}
One can also carry out some variant of the approach presented here, in order to obtain positivity of some positive combination of Schur forms which is not covered by what we do here. 

More specifically, one can show for instance the positivity of Schur forms such as $ c_2 $ in every rank and $ c_1 c_2 - c_3 $ for a rank $3$ Griffiths (semi)positive holomorphic Hermitian vector bundle. See \cite{Fag21} for more details.
\end{rem}
		
We obtain thus a partial affirmative answer to Griffiths' question for the polynomials in the Chern forms of $(E,h)$ belonging to the positive convex sub-cone $\mathcal F(r)\subset\Pi(r)$ spanned by \emph{all possible wedge products of all possible push-forwards} $\pis c_1(Q_{\rho}^\mathbf a)^{d_{\rho} + k}$, for $k=0,\dots,n$, as the weights $\mathbf a\in\mathbb N^r$ vary in the appropriate range prescribed by Conditions \eqref{eq:condition on a} and \eqref{eq:star}. This sub-cone contains in particular the signed Segre forms, which arise in the case of projectivized bundle.

\begin{rem}
Even if it is possible to obtain every Schur form as a push-forward, as observed in Remark \ref{rem:schurfunctions}, unfortunately this is not enough to get here that every Schur form of a Griffiths semipositive vector bundle is positive. This is because the curvature computations for the tautological bundles over the compete flag bundle do not permit to conclude that the relevant monomials whose push-forward give the desired Schur forms are positive.
\end{rem}

\medskip

The next examples are intended to give some flavor of which kind of new positive forms, in particular besides signed Segre forms, we are able to obtain with our methods.

	\subsection{Examples}\label{subsec:examples}
	In this section we give several examples of differential forms whose (strong) positivity is due to Theorem \ref{thm: positivity of push-forwards by flag bundles} and was not previously known in general. The explicit forms of some of them are obtained by implementing Formula \eqref{eq: darondeau-pragacz classes coefficient} in PARI/GP.
	
	\medskip
	
As we have seen, for $(E,h)$ a rank $r$ Griffiths (semi)positive vector bundle, the forms $\pi_*c_1(Q_{\rho}^\mathbf a,h)^{d_{\rho} + k}$ are strongly positive. We want to highlight here some among them that cannot be shown to be positive only using results in the literature preceding the present work, at the best of our knowledge.

To this aim, observe that we already knew that the signed Segre forms $(-1)^k s_k(E,h)$ are (strongly) positive for Griffiths (semi)positive vector bundles thanks to \cite[Theorem 1.1]{Gul12} (even though the strong positivity was not explicitly observed there). Also, as already noted, the product of positive forms, all of them strongly positive (resp. all except possibly one) is strongly positive (resp. positive). 

This understood, in order to check for which forms we get new information about their positivity, we shall express $\pi_*c_1(Q_{\rho}^\mathbf a,h)^{d_{\rho} + k}$ as a polynomial in the (signed) Segre forms of $E$.
	
In what follows, in order to simplify the notation, we denote by $ c_1,\ldots,c_r $ the Chern forms of $ (E,h) $ and by $ s_1,\ldots,s_n $ the Segre forms of $ (E,h) $.
The symbol $ S_{\sigma} $ stands for the Schur form of $ (E,h) $ associated to the partition $ \sigma $.
Moreover, we omit the symbol $ \wedge $ for the wedge product of forms.
	
	\subsubsection{Push-forwards from Grassmannian bundles}
	Denote by $ \rho $ the sequence $ (0,r-d,r) $. Then, $ \FF(E) $ is the Grassmannian bundle $ \GG_{r-d}(E) $ of $ (r-d) $-planes in $ E $.
	Let $ \pi\colon \GG_{r-d}(E)\to X$ be the projection, and denote by $ Q $ the universal quotient bundle of rank $ d $ on $ \GG_{r-d}(E) $ equipped with the quotient metric. In our notation, the class $ c_1(Q) $ equals $ c_1(Q_{\rho,1}) $.
	
	Therefore, for $ N \ge d(r-d) $ the metric counterpart of the Darondeau--Pragacz push-forward formula reads
	\begin{equation}\label{eq: push plucker class grassmannian}
	\pis c_1(Q,h)^{N} = \sum_{|\lambda| = N - d(r-d)} f^{\lambda + \varepsilon} \det\left( (-1)^{\lambda_i +j - i} s_{\lambda_i +j - i} \right)_{1 \le i,j \le d}
	\end{equation}
	where $ \lambda = (\lambda_1,\ldots,\lambda_d) $ is a partition and $ |\lambda| $ is its total weight, $ \varepsilon $ stands for the $ d $-uple $ (r-d)^d = (r-d,\ldots,r-d) $ and $ f^{\lambda + \varepsilon} $ is the number of standard Young tableaux with shape $ \lambda + \varepsilon $ (we have used here the more explicit version computed in the particular case of Grassmannian bundles in \cite[Theorem 0.1]{KT15}).
	Now, as explained for instance in \cite{KT15}, we have that
	\[
	f^{\lambda + \varepsilon} = \frac{N! \prod_{1\le i < j \le d}(\lambda_i - \lambda_j -i+j)}{\prod_{1\le i \le d}(r+\lambda_i -i)!} .
	\]
	
	Note that when $ d=1 $ the bundle $ \GG_{r-1}(E) $ can be identified with $ \PP(E^{\vee}) $, consequently $ Q \cong \OO_{E^{\vee}}(1) $, and Formula \eqref{eq: push plucker class grassmannian} becomes
	\[
	\pis c_1(\OO_{E^{\vee}}(1),h)^{N} = (-1)^{N-r+1} s_{N-r+1},
	\]
	which is the push-forward formula by \cite{Mou04,Gul12,Div16} giving the positivity of signed Segre forms.
	
	\medskip
	
	It is noteworthy (see Example \ref{ex:toycase}) to observe that already for the projectivized bundle of lines $ \PP(E) $ corresponding to the partition $ (0,1,r) $, if we push-forward powers of $ c_1(Q,h) $, where $ Q = \pi^{*}E / \OO_{E}(-1) $, we are now able to get forms whose positivity was not previously known.
	
	In rank $ 3 $ (if $ r=2 $ we have that $ \PP(E) \cong \PP(E^{\vee}) $ and there is nothing more to add) we see for instance, by using Formula \eqref{eq: push plucker class grassmannian}, that $ \pis c_1(Q,h)^5 $ equals the form
	$$
	4 c_1^3 - 3 c_1c_2 - c_3 = s_3 - 5 s_1s_2 ,
	$$
	and the positivity of $ s_3 - 5 s_1s_2 $ was not previously known given that $ s_3 $ is negative.
	
	Analogously, for the same reasons, if $ r=4 $, the positivity of
	\[
	\pis c_1(Q,h)^6 = 10 c_1^3 - 4c_1c_2 - c_3 = s_3 -6 s_1s_2 - 5 s_1^3
	\]
	was not previously known.
	
	\medskip
	
	The simplest example of a Grassmannian bundle which is not a projectivized bundle is $ \GG_{2}(E) $ for $ E $ of rank $ 4 $. Also in this case we get something new, as follows.
	By Formula \eqref{eq: push plucker class grassmannian}, the push-forward \textsl{via} $ \pi $ of $ c_1(Q,h)^N $ is given by
	$$
	\begin{cases}
		 2 & \textrm{for $ N=4 $,} \\ 
		5c_1&  \textrm{for $ N=5 $,} \\
		9c_1^2 - 4c_2 &  \textrm{for $ N=6 $,} \\
		14(c_1^3 -  c_1c_2) &  \textrm{for $ N=7 $,} \\
		2(10 c_1^4 - 16 c_1^2 c_2 - c_1 c_3 + 3 c_2^2 + 4 c_4) &  \textrm{for $ N=8 $.}
	\end{cases}
	$$
	When rewritten in terms of signed Segre forms, we obtain
	$$
	\begin{cases}
		 2 & \textrm{for $ N=4 $,} \\ 
		5(-s_1)&  \textrm{for $ N=5 $,} \\
		5s_1^2+4s_2 &  \textrm{for $ N=6 $,} \\
		14(-s_1)s_2&  \textrm{for $ N=7 $,} \\
		2(7 s_1s_3 + 7 s_2^2 - 4 s_4) &  \textrm{for $ N=8 $,}
	\end{cases}
	$$
	so that the positivity of the last form could not be previously deduced, since $ -s_4 $ is negative.
	
	\subsubsection{Push-forwards from complete flag bundles: the case of rank $3$}
	The general formul\ae{} for the push-forwards of $ c_1(Q^{(a,b,c)},h)^{3+k} $ in terms of the Schur forms, up to degree $ 3 $, are:
$$
\begin{aligned}
	&k=0\quad\leadsto\quad  3 (a^2b - ab^2 - a^2c + ac^2 + b^2 c - bc^2),\\
	&k=1\quad\leadsto\quad	  4(a^3b -ab^3 - a^3c + ac^3 + cb^3 - bc^3) S_{(1)}, \\
	&k=2\quad\leadsto\quad 10(a^3b^2 - a^2b^3 - a^3c^2 + a^2c^3 + b^3c^2 - b^2c^3) S_{(2,0)} \\
		&  \qquad\qquad\qquad\qquad+ 5(a^4b - ab^4 - a^4c + ac^4 + b^4c - bc^4) S_{(1,1)}, \\
	&k=3\quad\leadsto\quad	60(a^3b^2c - a^2b^3c - a^3bc^2 + a^2bc^3 + ab^3c^2 -  ab^2c^3)S_{(3,0,0)} \\
	&  \qquad\qquad\qquad\qquad + 15(a^4b^2 - a^2b^4 - a^4c^2 + a^2c^4 + b^4c^2 - b^2c^4) S_{(2,1,0)} \\
		 &\qquad\qquad\qquad\qquad\qquad +6 (a^5b - ab^5 - a^5c + ac^5 + b^5c - bc^5) S_{(1,1,1)} .
\end{aligned}
$$
	
	Clearly, our method produces positive forms until $ k $ reaches $ n $, but already from these first cases we see how the complexity rapidly increases.
	Thanks to our Theorem \ref{thm: positivity of push-forwards by flag bundles} we can say that all of the above listed forms belong to the positive cone for every $ a \ge b \ge c \ge 0 $.
	
	For a concrete example, if $ (a,b,c) = (3,2,0) $ we obtain in terms of Segre forms:
	\begin{align*}
		\pis c_1(Q^{(3,2,0)},h)^{6} &= 2700 \,S_{(2,1,0)} + 2340\, S_{(1,1,1)} \\
		&= 180 (- 15 s_1s_2 + 2 s_3),
	\end{align*}
	and the positivity of this form was not already known since $ s_3 < 0 $.
	
	On the negative side, note that, for instance, the positivity of 
	$$ 
	\pis c_1(Q^{(2,1,0)},h)^{6} = -180 s_1 s_2 
	$$ 
	was instead previously known since $ - s_1 s_2 $ is the wedge product of strongly positive forms.
	
	\subsubsection{Push-forwards from complete flag bundles: the case of rank $4$}
	Here, we prefer to emphasize the different behavior in two special cases instead of giving the general formul\ae{} for $ a\ge b\ge c\ge d\ge 0 $.
	
	For $ (a,b,c,d) = (3,2,1,0) $, we have
	\begin{align*}
	&\pis c_1(Q^{(3,2,1,0)},h)^{9} = 90720 ( - s_1^3 - 2 s_1 s_2), \\
	&\pis c_1(Q^{(3,2,1,0)},h)^{10} = 5040 (216 s_1^2 s_2 + 7 s_1 s_3 + 39 s_2^2 - 4 s_4),
	\end{align*}
	and note that the positivity of the last form was not previously known since $ -s_4 $ is negative, while we already knew that $ - s_1^3 - 2 s_1 s_2 $ is positive.
	
	If we want an example of a push-forward that gives not previously known positive forms in degree $ 3 $ and $ 4 $, consider the case $ (a,b,c,d) = (4,3,2,0) $ where we have
	\begin{align*}
	&\pis c_1(Q^{(4,3,2,0)},h)^{9} = 181440 (-8 s_1^3 - 12 s_1 s_2 + s_3), \\
	&\pis c_1(Q^{(4,3,2,0)},h)^{10} = 40320 ( 648 s_1^2 s_2 - 124 s_1 s_3 + 42 s_2^2 + 13 s_4 )
	\end{align*}
	and the positivity of both of these forms was not already known again because $ s_3 $ is negative.

\bibliography{bibliography}{}

\providecommand{\bysame}{\leavevmode\hbox to3em{\hrulefill}\thinspace}
\providecommand{\MR}{\relax\ifhmode\unskip\space\fi MR }
\providecommand{\MRhref}[2]{%
  \href{http://www.ams.org/mathscinet-getitem?mr=#1}{#2}
}
\providecommand{\href}[2]{#2}
\begin{thebibliography}{Dem88b}

\bibitem[Dam73]{Dam73}
James Damon, \emph{The {G}ysin homomorphism for flag bundles}, Amer. J. Math.
  \textbf{95} (1973), 643--659. \MR{348760}

\bibitem[Dem88a]{Dem88}
Jean-Pierre Demailly, \emph{Vanishing theorems for tensor powers of a positive
  vector bundle}, Geometry and analysis on manifolds ({K}atata/{K}yoto, 1987),
  Lecture Notes in Math., vol. 1339, Springer, Berlin, 1988, pp.~86--105.
  \MR{961475}

\bibitem[Dem88b]{Dem88amp}
\bysame, \emph{Vanishing theorems for tensor powers of an ample vector bundle},
  Invent. Math. \textbf{91} (1988), no.~1, 203--220. \MR{918242}

\bibitem[Dem12]{Dem01}
\bysame, \emph{{C}omplex {A}nalytic and {D}ifferential {G}eometry}, Available
  at
  \texttt{https://www-fourier.ujf-grenoble.fr/\~demailly/manuscripts/agbook.pdf},
  Version of 2012.

\bibitem[Div16]{Div16}
Simone Diverio, \emph{Segre forms and {K}obayashi-{L}\"{u}bke inequality},
  Math. Z. \textbf{283} (2016), no.~3-4, 1033--1047. \MR{3519994}

\bibitem[DP17]{DP17}
Lionel Darondeau and Piotr Pragacz, \emph{Universal {G}ysin formulas for flag
  bundles}, Internat. J. Math. \textbf{28} (2017), no.~11, 1750077, 23 pp.
  \MR{3714353}

\bibitem[DPS94]{DPS94}
Jean-Pierre Demailly, Thomas Peternell, and Michael Schneider, \emph{Compact
  complex manifolds with numerically effective tangent bundles}, J. Algebraic
  Geom. \textbf{3} (1994), no.~2, 295--345. \MR{1257325}

\bibitem[Fag22]{Fag21}
Filippo Fagioli, \emph{A note on {G}riffiths' conjecture about the positivity
  of {C}hern-{W}eil forms}, Differential Geom. Appl. \textbf{81} (2022), Paper
  No. 101848, 15. \MR{4375642}

\bibitem[Fin20]{Fin20}
Siarhei Finski, \emph{{On characteristic forms of positive vector bundles,
  mixed discriminants and pushforward identities}}, arXiv e-prints,
  arXiv:2009.13107, 2020.

\bibitem[FL83]{FL83}
William Fulton and Robert Lazarsfeld, \emph{Positive polynomials for ample
  vector bundles}, Ann. of Math. (2) \textbf{118} (1983), no.~1, 35--60.
  \MR{707160}

\bibitem[GH78]{GH78}
Phillip Griffiths and Joseph Harris, \emph{Principles of algebraic geometry},
  Wiley-Interscience [John Wiley \& Sons], New York, 1978, Pure and Applied
  Mathematics. \MR{507725}

\bibitem[Gri69]{Gri69}
Phillip~A. Griffiths, \emph{Hermitian differential geometry, {C}hern classes,
  and positive vector bundles}, Global {A}nalysis ({P}apers in {H}onor of {K}.
  {K}odaira), Univ. Tokyo Press, Tokyo, 1969, pp.~185--251. \MR{0258070}

\bibitem[Gul12]{Gul12}
Dincer Guler, \emph{On {S}egre forms of positive vector bundles}, Canad. Math.
  Bull. \textbf{55} (2012), no.~1, 108--113. \MR{2932990}

\bibitem[Ilo78]{Ilo78}
Samuel~A. Ilori, \emph{A generalization of the {G}ysin homomorphism for flag
  bundles}, Amer. J. Math. \textbf{100} (1978), no.~3, 621--630. \MR{501229}

\bibitem[KT15]{KT15}
Hajime Kaji and Tomohide Terasoma, \emph{Degree formula for {G}rassmann
  bundles}, J. Pure Appl. Algebra \textbf{219} (2015), no.~12, 5426--5428.
  \MR{3390031}

\bibitem[Li20]{Li20}
Ping Li, \emph{{Nonnegative Hermitian vector bundles and Chern numbers}},
  Mathematische Annalen (2020).

\bibitem[Mou04]{Mou04}
Christophe Mourougane, \emph{Computations of {B}ott-{C}hern classes on {${\Bbb
  P}(E)$}}, Duke Math. J. \textbf{124} (2004), no.~2, 389--420. \MR{2079253}

\bibitem[Pin18]{Pin18}
Vamsi~Pritham Pingali, \emph{Representability of {C}hern-{W}eil forms}, Math.
  Z. \textbf{288} (2018), no.~1-2, 629--641. \MR{3774428}

\bibitem[Xia20]{Xia20}
Jian Xiao, \emph{{On the positivity of high-degree Schur classes of an ample
  vector bundle}}, arXiv e-prints, arXiv:2007.12425, 2020.

\end{thebibliography}

\end{document}